\documentclass[a4paper,11pt]{article}
\usepackage[french,english]{babel}

\usepackage[utf8]{inputenc}

\usepackage{amssymb,amscd,amsmath, amsthm, amsfonts,graphicx}
\usepackage[pdftex,colorlinks, citecolor=blue, linkcolor=red]{hyperref}
\hypersetup{hypertexnames=false}
\usepackage[capitalise]{cleveref}
\usepackage[utf8]{inputenc}
\numberwithin{equation}{section}
\newtheorem{theorem}{Theorem}
\newtheorem{lemma}[theorem]{Lemma}

\newtheorem{definition}[theorem]{Definition}

\newtheorem{remark}[theorem]{Remark}
\numberwithin{theorem}{section}
\crefformat{equation}{(#2#1#3)}

\begin{document}

\title{Joint spectrum and large deviation principle for random matrix products}
\author{Cagri SERT}
\date{}
\maketitle



\begin{abstract}
The aim of this note is to announce some results about the probabilistic and deterministic asymptotic properties of linear groups. The first one is the analogue, for norms of random matrix products, of the classical theorem of Cram\'{e}r on large deviation principles (LDP) for sums of iid real random variables. In the second result, we introduce a limit set describing the asymptotic shape of the powers $S^{n}=\{g_{1}. \ldots.g_{n} \, | \, g_{i} \in S\}$ of a subset $S$ of a semisimple linear Lie group $G$ (e.g. SL$(d,\mathbb{R})$). This limit set has applications, among others, in the study of large deviations.
\end{abstract}

\section{Large deviation principle for random matrix products}
\subsection{Introduction}
Let $G$ be a connected semisimple linear real algebraic group 
(e.g. SL$(d,\mathbb{R})$). A random walk on $G$ is a random 
process $Y_{n}=X_{n}.\ldots.X_{1}$ where $X_{i}$'s are 
independent and identically distributed (iid) $G$-valued 
random variables. Starting from Bellman \cite{Bellman}, 
Furstenberg-Kesten \cite{Furstenberg.Kesten} and Furstenberg \cite{Furstenberg.non.commuting}, an important aim in the 
study of these non-commutative random walks was to establish 
the analogues of the classical limit theorems existing for 
the iid real random variables. More precisely, one is 
interested in studying the probabilistic limiting behavior of 
the (logarithms of the) norms $\log||Y_{n}||$ of random 
matrix products. In fact, we will consider the slightly more 
general multi-norm given by the Cartan projection that we 
briefly describe before going on: let $\mathfrak{g}$ be the 
Lie algebra of the group $G$, $\mathfrak{a}$ a Cartan 
subalgebra of $\mathfrak{g}$ and $\mathfrak{a}^{+}$ a chosen 
Weyl chamber. Let $K$ be a maximal compact subgroup of $G$ 
for which we have the Cartan decomposition $G=K\exp(\mathfrak{a}^{+})K$. Then the map $\kappa: G \longrightarrow \mathfrak{a}^{+}$ is well-defined by the 
following equality and is called the Cartan projection: for 
all $g \in G$, $g=k\exp(\kappa(g))u$ for some $k,u \in K$. As 
an example, in the case of $G=\text{SL}(d,\mathbb{R})$, the Cartan 
projection $\kappa(g)$ of a matrix $g \in \text{SL}(d,\mathbb{R})$ 
writes as $\kappa(g)=(\log ||g||, \log \frac{||\wedge^{2}g||}{||g||}, \ldots, \log \frac{||\wedge^{d}g||}{||\wedge^{d-1}g||})$, where $\wedge^{k}\mathbb{R}^{d}$'s are considered 
with their canonical Euclidean structures and $||.||$'s 
denote the associated operator norms. The components of $\kappa(g)$ are the logarithms of the singular values of $g$.

The first limit theorem that was proven for random matrix 
products was the analogue of the law of large numbers. 
Stating it in our setting, Furstenberg-Kesten's result   \cite{Furstenberg.Kesten} reads: if $\mu$ is a probability 
measure on $G$ with a finite first moment (i.e. $\int ||\kappa(g)|| \mu(dg) < \infty$ for some norm $||.||$ on $\mathfrak{a}$), 
then the $\mu$-random walk $Y_{n}=X_{n}.\ldots.X_{1}$ (i.e. $X_{i}$'s are iid of law $\mu$) satisfies $$\frac{1}{n} \kappa(Y_{n}) \overset{a.s.}{\underset{n \rightarrow \infty}{\longrightarrow}} \vec{\lambda}_{\mu} \in \mathfrak{a}$$ where $\vec{\lambda}_{\mu}$ can be defined by 
this and called the Lyapunov vector of $\mu$. Nowadays, this 
result is a rather straightforward corollary of Kingman's 
subadditive ergodic theorem. A second important limit theorem 
that was established in increasing generality by Tutubalin \cite{Tutubalin}, Le Page \cite{LePage}, Goldsheid-Guivarc'h \cite{Guivarch.Goldsheid}, and Benoist-Quint \cite{Benoist.Quint.CLT} \cite{BQpoly} is the central limit 
theorem (CLT). Benoist-Quint's CLT reads: if $\mu$ is a 
probability measure on $G$ of finite second order moment and 
such that the support of $\mu$ generates a Zariski-dense 
semigroup in $G$, then $\frac{1}{\sqrt{n}}(\kappa(Y_{n})-n\vec{\lambda}_{\mu})$ converges in distribution to a non-
degenerate Gaussian law on $\mathfrak{a}$. A feature of this 
result is the Zariski density assumption which also appears in our result below. We note that the fact that the support $S$ of 
the probability measure $\mu$ generates a Zariski-dense 
semigroup can be read as: any polynomial that vanishes on $\cup_{n \geq 1} S^{n}$ also vanishes on $G$ (we recall that 
in $\mathbb{R}$, a subset is Zariski dense if and only if it 
is infinite). Some other limit theorems whose analogues have 
been obtained are the law of iterated logarithm and local 
limit theorems for which we refer the reader to the nice 
books of Bougerol-Lacroix \cite{Bougerol.Lacroix} and more 
recently Benoist-Quint \cite{BQpoly}.

An essential and, up to our work, a rather incomplete aspect of these non-commutative limit theorems is concerned with large deviations. The main result in this direction is 
that of Le Page \cite{LePage}, (see also Bougerol \cite{Bougerol.Lacroix}) and its extension by Benoist-Quint \cite{BQpoly}, stating the exponential decay of large 
deviation probabilities off the Lyapunov vector. Before 
stating this result, recall that a probability measure $\mu$ 
on $G$ is said to have a finite exponential moment, if there 
exists $\alpha>1$ such that we have $\int \alpha^{||\kappa(g)||} \mu(dg)<\infty$. We have 
\begin{theorem}[Le Page \cite{LePage}, Benoist-Quint \cite{BQpoly}] \label{LePage}
Let $G$ be as before, $\mu$ be a probability measure of 
finite exponential moment on $G$ whose support generates a 
Zariski-dense semigroup in $G$. Let $Y_{n}$ denote the $n^{th}$ step of the $\mu$-random walk on $G$. Then, for all $\epsilon>0$, we have $\limsup_{n \rightarrow \infty} \frac{1}{n} \log \mathbb{P}(||\frac{\kappa(Y_{n})}{n}-\vec{\lambda}_{\mu}||>\epsilon)<0$.
\end{theorem}
\subsection{Statement of main result}
In our first main result, under the usual Zariski density 
assumption, we prove the matrix analogue of a 
classical theorem (see below) about large deviations for iid real 
random variables. Let $X$ be a topological space and $\mathcal{F}$ be a $\sigma$-
algebra on $X$.
\begin{definition} A sequence $Z_{n}$ of $X$-valued random 
variables is said to satisfy a large deviation principle 
(LDP) with rate function $I:X \longrightarrow [0, \infty]$, 
if for every measurable subset $R$ of $X$, we have $$
\underset{x \in \overset{\circ}{R}}{-\inf I(x)} \leq \underset{n \rightarrow \infty}{\liminf} \frac{1}{n}\log \mathbb{P}(Z_{n} \in R) \leq \underset{n \rightarrow \infty}{\limsup} \frac{1}{n}\log \mathbb{P}(Z_{n} \in R) \leq \underset{x \in \overline{R}}{-\inf I(x)} $$ 
\end{definition} \noindent where, $\overset{\circ}{R}$ denotes the interior and $\overline{R}$ the closure of $R$. 

With this definition, the classical Cram\'{e}r-Chernoff 
theorem says that the sequence of averages $Y_{n}=\frac{1}{n}\sum_{i=1}^{n}X_{i}$ of real iid random variables of finite 
exponential moment satisfies an LDP with a proper convex rate 
function $I$, given by the convex conjugate (Fenchel-Legendre 
transform) of the logarithmic moment generating function of $X_{i}$'s.  Our first main result reads
\begin{theorem} \label{kanitsav1}
Let $G$ be a connected semisimple linear real algebraic group 
and $\mu$ be a probability measure of finite exponential 
moment on $G$, whose support generates a Zariski dense sub-
semigroup of $G$. Then, the sequence random variables $\frac{1}{n}\kappa(Y_{n})$ satisfies an LDP with a proper 
convex rate function $I:\mathfrak{a} \longrightarrow [0,\infty]$ assuming a unique zero on the Lyapunov vector $\vec{\lambda}_{\mu}$ of $\mu$.
\end{theorem}
\begin{remark}\label{remark1}
1. Without any moment assumptions on $\mu$, we also obtain a weaker result which is an analogue of a result of Bahadur \cite{Bahadur} for iid real random variables.\\[4pt]
2. Under a stronger exponential moment condition (sometimes 
called finite super-exponential moment), by exploiting 
convexity of $I$, we are able to identify the rate function $I$ with the convex conjugate of a limiting logarithmic 
moment generating function of the random variables $\frac{1}{n}\kappa(Y_{n})$.\\[4pt]
3.  It follows by convexity of $I$ that the effective support $D_{I}:=\{x \in \mathfrak{a} \; | \; I(x) < \infty \}$ of the 
rate function $I$ is a convex subset of $\mathfrak{a}^{+}$. 
We also show that this set $D_{I}$ depends only on the 
support $S$ of $\mu$ (i.e. not on the particular mass 
distribution $\mu$ on $S$) and we identify $D_{I}$ with a 
Hausdorff limiting set of a deterministic construction, 
namely the joint spectrum $J(S)$ of $S$ (see Theorem \ref{kanitsav3} below).\\[4pt]
4.  We conjecture that a similar LDP holds for the Jordan projection $\lambda: G \to \mathfrak{a}^{+}$ in place of $\kappa$ (see Def. in paragraph 2.1).
\end{remark}

\subsection{The idea of the proof}

The fact that $I(\vec{\lambda}_{\mu})=0$ follows from the definition of $\vec{\lambda}_{\mu}$ and $I$. That this zero is unique is merely the 
translation of the deep result of Le Page (Theorem \ref{LePage}). 
That $I$ is proper is a rather straightforward consequence of 
Chernoff's estimates and that one can identify $I$ as a 
Fenchel-Legendre transform (2. of Remark \ref{remark1}) follows 
from its convexity using standard techniques, namely the 
Varadhan's integral lemma (see Theorem 4.3.1. in \cite{Dembo.Zeitouni}) and elementary properties of convex 
conjugation of functions. Moreover, the convexity of the rate 
function is proven using ideas similar to those used to prove its existence. Therefore, in this note, we will focus on the proof of the existence of an LDP. For simplicity, we shall assume that the measure $\mu$ is compactly 
supported. We make use of the following general fact:
\begin{theorem}[see Theorem 4.1.11 in \cite{Dembo.Zeitouni}] \label{DZ.ktsav}
Let $X$ be a topological space endowed with its Borel $\sigma$-algebra $\beta_{X}$, and $Z_{n}$ be a sequence of $X$-
valued random variables taking values in a compact subset of $X$. Let $\mathcal{A}$ be a base of open sets for the topology of $X$. For 
each $x \in X$, define:
$$
I_{li}(x):= \underset{\underset{x \in A}{A \in \mathcal{A}}} {\sup} - \underset{n \rightarrow \infty}{\liminf} \frac{1}{n}\log \mathbb{P}(Z_{n} \in A) \qquad \text{and} \qquad I_{ls}(x):= \underset{\underset{x \in A}{A \in \mathcal{A}}}{\sup} - \underset{n \rightarrow \infty}   {\limsup} \frac{1}{n}\log \mathbb{P}(Z_{n} \in A) 
$$
Suppose that for all $x \in X$, we have $I_{li}(x)=I_{ls}(x)$. Then, the sequence $Z_{n}$ satisfies an LDP with rate 
function $I$ given by $I(x):=I_{li}(x)=I_{ls}(x)$.
\end{theorem}
In view of this theorem, to prove the existence of the LDP in Theorem \ref{kanitsav1}, we have to show that the equality $I_{li}=I_{ls}$ is satisfied. A first remark is if $\kappa$ was an additive mapping (i.e. $\kappa(gh)=\kappa(g)+\kappa(h)$), this would follow rather easily from the independence of random walk increments and uniform continuity of $\kappa$. A further important remark is that in fact a weaker form of additivity (i.e. $||\kappa(gh)-\kappa(g)-\kappa(h)||$ is uniformly bounded for all $g,h \in \text{supp}(\mu)$) is sufficient to insure this equality. A key result of Benoist \cite{Benoist1}, \cite{Benoist2} shows that this weak form of additivity  is satisfied in a given $(r,\epsilon)$-Schottky semigroup. This already finishes the proof in the case when $\mu$ is supported on such a semigroup. For the general case, we need an argument showing that we can restrict the random walk on Schottky semigroups with no loss in the exponential rate of decay probabilities (as in Lemma \ref{AMS.dispersion}). This is 
done by using a result of Abels-Margulis-Soifer \cite{AMS} about proximal elements in Zariski dense semigroups (which in 
turn uses a result of Benoist-Labourie \cite{Benoist.Labourie} and Prasad \cite{Prasad}) together 
with the uniform continuity of the Cartan projection and then 
a simple partitioning and pigeonhole argument. A key step in this proof is the following Lemma \ref{AMS.dispersion}. 

Recall that an element $g$ in $\text{PGL}(d,\mathbb{R})$ is called $\epsilon$-proximal if it has a unique eigenvalue of maximal modulus and the ratio of its first two singular values is at least $\frac{1}{\epsilon}$. It is called $(r,\epsilon)$-proximal if additionally the top eigenvector is at least $r$ away from projective hyperplane spanned by the other generalized eigenspaces. An element $g$ in $G$ is called $(r,\epsilon)$-loxodromic if it is $(r,\epsilon)$-proximal in each of the rk$(G)$ fundamental proximal representations of G. 

Abels-Margulis-Soifer show that for a Zariski 
dense semigroup $\Gamma$ in $G$, there exists $r>0$ such that 
for every $\epsilon>0$, one can find a \textit{finite} subset $F \subset \Gamma$ with the property that for all $\gamma \in \Gamma$, there exists $f \in F$ such that $\gamma.f$ is $(r,\epsilon)$-loxodromic. We denote by $\Gamma$ the semigroup 
generated by the support of $\mu$.
\begin{lemma}\label{AMS.dispersion} Let $0<\epsilon<r=r(\Gamma)$. There exist a compact set $C=C(\Gamma,\epsilon) \subset \mathfrak{a}$, a natural number $i_{0}=i_{0}(\epsilon,\Gamma, \mu)$, and a constant $d_{1}=d_{1}(\epsilon,\Gamma,\mu)>0$ such that for all $n_{0} \in \mathbb{N}$ and subset $R \subset \mathfrak{a}^{+}$, there exists a 
natural number $n_{1} \geq n_{0}$ with $n_{1}-n_{0}\leq i_{0}$ such that we have
$$
\mathbb{P}( \kappa (Y_{n_{1}}) \in R + C \; \text{and} \; Y_{n_{1}} \; \text{is} \; (r,\epsilon)\text{-loxodromic}) \geq d_{1}.\mathbb{P}(\kappa(Y_{n_{0}}) \in R)
$$
\end{lemma}

The next step in the proof consists in observing that one can further restrict the random 
walk to an $(r,\epsilon)$-Schottky semigroup, again losing only 
a uniformly bounded proportion of the probability (as in Lemma  \ref{AMS.dispersion}). By doing so, we reduce the situation 
to a random walk on a semigroup on which the Cartan 
projection $\kappa(.)$ is almost additive and hence we 
can conclude as we sketched in the beginning of the argument.

\begin{remark}
In a further work in preparation \cite{Sert.hyperbolic.LDP}, using 
essentially the same method as in the proof of Theorem \ref{kanitsav1}, we show that the LDP holds for the 
average word length of random walks on Gromov hyperbolic 
groups.
\end{remark}

\section{Joint spectrum}
\subsection{Introduction}
In this second part of this note, we define the notion of \emph{joint spectrum} of a bounded subset of a semisimple Lie group. We then relate this notion to the effective support $D_{I}=\{x \in \mathfrak{a}\, | \, I(x)<\infty \}$ of the rate function from Theorem \ref{kanitsav1}. Recall the definition of the Jordan projection $\lambda: G \longrightarrow \mathfrak{a}^{+}$: if $g=g_{e}g_{h}g_{u}$ is the Jordan decomposition of $g$ with $g_{e}$ elliptic, $g_{h}$ hyperbolic and $g_{u}$ is unipotent, then $\lambda(g)$ is defined as $\kappa(g_{h})$. Let $S$ be a bounded subset of $G$ and $S^{n}=\{g_{1}.\ldots.g_{n} \, | \, g_{i} \in S\}$ denote its $n^{th}$ power. We are interested in the following questions: do the sequences of bounded subsets of $\mathfrak{a}^{+}$, $\frac{1}{n}\kappa(S^{n})$ and $\frac{1}{n}\lambda(S^{n})$ have a limit in the Hausdorff topology ? If yes, are the limits the same and can one describe these limit sets ?
 
\subsection{Statement of the main results}
Regarding the above questions, we show the following: 

\begin{theorem}\label{kanitsav2}
Let $G$ be a connected semisimple linear real algebraic group 
and $S$ a bounded subset of $G$ generating a Zariski dense 
sub-semigroup. Then,\\[4pt]
1. The following Hausdorff limits exist, and we have the 
equality: $$ \lim_{n \rightarrow \infty}\frac{1}{n}\kappa(S^{n})=\lim_{n \rightarrow \infty} \frac{1}{n}\lambda(S^{n})$$ This common limit will be denoted as $J(S)$ 
and called the \emph{joint spectrum} of $S$. \\[4pt]
2. $J(S)$ is a compact convex subset of $\mathfrak{a}^{+}$ with non-empty interior.
\end{theorem}

\begin{remark} \label{remark2}
1. It is not hard to see that if two subsets $S$ and $S'$ of $G$ generate the same Zariski dense semigroup $\Gamma$ in $G$, then the projective images of $J(S)$ and $J(S')$ are the 
same. Therefore the corresponding cone in $\mathfrak{a}^{+}$ only depends on $\Gamma$. It turns out that this 
is precisely the Benoist limit cone of $\Gamma$ \cite{Benoist2}. \\[4pt]
2. For a Banach algebra $(\mathcal{B}, ||.||)$ and a bounded 
subset $S$ of $\mathcal{B}$, denote the joint spectral radius 
of $S$ by $r(S):=\lim_{n \rightarrow \infty}\frac{1}{n}\log(\sup_{x \in S^{n}}||x||)$. When $\mathcal{B}$ is finite 
dimensional, this does not depend on the particular norm on $\mathcal{B}$. Let now $\rho:G \rightarrow \text{SL}(V)$ be an 
irreducible rational representation of $G$ and let $\chi_{\rho} \in \mathfrak{a}^{\ast}$ denote the highest 
weight of $\rho$. Then we have the equality $\sup_{z \in J(S)}\chi_{\rho}(z)=r(\rho(S))$. This allows us to derive a multi-dimensional generalization of 
Berger-Wang identity \cite{Berger.Wang}.
\end{remark}

We now relate the joint spectrum of 
the support of a probability measure $\mu$ with the effective 
support $D_{I}$ of the rate function $I$ for the LDP of $\frac{1}{n} \kappa(Y_{n})$, where $Y_{n}$ denotes as usual 
the $\mu$-random walk on $G$. 

\begin{theorem} \label{kanitsav3}
Let $G$ and $\mu$ be as in Theorem \ref{kanitsav1}, and suppose 
moreover that the support $S$ of $\mu$ is bounded. Let $I$ be 
the rate function given by Theorem \ref{kanitsav1} and $D_{I}=\{x \in \mathfrak{a} \; | \; I(x)<\infty\}$ be the effective 
support of $I$. We then have \\[1pt]
1. $\overline{D}_{I}=J(S)$ and $\overset{\circ}{D_{I}}=\overset{\circ}{J(S)}$. If $S$ is moreover finite, then $D_{I}=J(S)$.\\[3pt]
2. The Lyapunov vector $\vec{\lambda}_{\mu} \in \mathfrak{a}^{+}$ of $\mu$ belongs to the interior of $J(S)$.
\end{theorem}

For the point 2. of this theorem, we note that for a probability measure as 
in Theorem \ref{kanitsav1}, the fact that $\vec{\lambda}_{\mu}$ 
belongs to the interior of $\mathfrak{a}^{+}$ was obtained by 
Guivarc'h-Raugi \cite{Guivarch.Raugi} and Goldsheid-Margulis \cite{Goldsheid.Margulis}. Our result gives 
a more precise location for $\vec{\lambda}_{\mu}$ in case $\mu$ is, moreover, boundedly supported. Our method extends also to the case where $\mu$ only has a finite exponential moment to show that $\vec{\lambda}_{\mu}$ belongs to the Benoist cone of the semigroup generated by the support of $\mu$.

The tools that go into the proof of Theorem \ref{kanitsav2} and Theorem \ref{kanitsav3} are mostly similar to those used in the 
proof of Theorem \ref{kanitsav1}. We use an additional tool to prove 
the fact that the joint spectrum is of non-empty interior: 
while this result can be directly deduced from the properties 
of Benoist cone proved in \cite{Benoist2}, we adopt an 
indirect approach and use the central limit theorem of 
Goldsheid-Guivarc'h \cite{Guivarch.Goldsheid} and Guivarc'h  \cite{Guivarch.Goldsheid2} that we combine with Abels-
Margulis-Soifer result \cite{AMS}  and Benoist's estimates \cite{Benoist1}. This in turn allows us to derive point 2. of the previous theorem and hence also part of the point 2. of Theorem \ref{kanitsav2}.

\begin{remark}\label{growth.ind.remark}
The notion of joint spectrum plays also an 
important role in another work in preparation \cite{Sert.growth.indicator},    
 where we study a new exponential counting function for 
a finite subset $S$ in a group $G$ as before.
\end{remark}

\section*{Acknowledgements} 
These results are part of author's doctoral thesis realized under 
the supervision of Emmanuel Breuillard in Universit\'{e} 
Paris-Sud during 2013-2016. The author would like to take the 
opportunity to thank him for 
numerous discussions and for suggesting the notion of joint 
spectrum. The author also thanks to WWU M\"{u}nster where part of this work was conducted.

\small \indent E-mail address: cagri.sert@math.ethz.ch \\
\small \indent Telephone number: +41 44 632 3737


\begin{thebibliography}{99}
\thispagestyle{empty}

\bibitem{AMS} Abels, H., Margulis, G.A. and Soifer, G.A., 1995. Semigroups containing proximal linear maps. Israel journal of mathematics, 91(1-3), pp.1-30.\\[-20pt]

\bibitem{Benoist1} Benoist, Y., 1996. Actions propres sur les espaces homog\`{e}nes r\'{e}ductifs. Annals of mathematics, pp.315-347.\\[-20pt]

\bibitem{Benoist2} Benoist, Y., 1997. Propri\'{e}t\'{e}s asymptotiques des groupes lin\'{e}aires. Geometric \& Functional Analysis GAFA, 7(1), pp.1-47.\\[-20pt] 

\bibitem{Benoist3} Benoist, Y., 2000, Propri\'{e}t\'{e}s asymptotiques des groupes lin\'{e}aires II. Advanced Studies Pure Math. 26 (2000) pp.33-48.\\[-20pt]

\bibitem{BQpoly} Benoist, Y. and Quint, J.F., 2013. Random walks on reductive groups. to appear.\\[-20pt]

\bibitem{Benoist.Quint.CLT} Benoist, Y. and Quint, J.F., 2016. Central limit theorem for linear groups. Ann. Probab, 44(2), pp.1308-1340.\\[-20pt]

\bibitem{Bahadur} Bahadur, R.R., 1971. Some limit theorems in statistics (No. 04; QA276. A1, B3 1971.). Philadelphia: Society for industrial and applied mathematics.\\[-20pt]

\bibitem{Bellman} Bellman, R., 1954. Limit theorems for non-commutative operations. I. Duke Mathematical Journal, 21(3), pp.491-500.\\[-20pt]

\bibitem{Benoist.Labourie}  Benoist, Y., and Labourie, F., 1993. Sur les diff\'{e}omorphismes d'Anosov affines \`{a} feuilletages stable et instable diff\'{e}rentiables. Inventiones mathematicae, 111(1), pp.285-308.\\[-20pt]

\bibitem{Berger.Wang} Berger, M.A. and Wang, Y., 1992. Bounded semigroups of matrices. Linear Algebra and its Applications, 166, pp.21-27.\\[-20pt]

\bibitem{Bougerol.Lacroix} Bougerol, P. and Lacroix J., 1985. Products of random matrices with applications to Schr\"{o}dinger operators. Progress in Probability and Statistics, 8. Birkh\"{a}user Boston, Inc., Boston, MA, xii+283 pp.\\[-20pt]

\bibitem{Dembo.Zeitouni} Dembo, A. and Zeitouni, O. Large deviations techniques and applications (Vol. 38). Springer Science \& Business Media. 2009\\[-20pt]

\bibitem{Furstenberg.non.commuting} Furstenberg, H. Noncommuting random products. Transactions of the American Mathematical Society, 1963, 108(3), pp.377-428.\\[-20pt]

\bibitem{Furstenberg.Kesten} Furstenberg, H. and Kesten, H. Products of random matrices. The Annals of Mathematical Statistics,  1960, 31(2), pp.457-469.\\[-20pt]

\bibitem{Goldsheid.Margulis} Goldsheid, I. and Margulis, G. Lyapunov Indices of a Product of Random Matrices, Russian Math. Surveys 44, 1989, pp.11-81.\\[-20pt]

\bibitem{Guivarch.Goldsheid} Goldsheid, I.Y. and Guivarc'h, Y. Zariski closure and the dimension of the Gaussian low of the product of random matrices. I. Probability theory and related fields, 1996, 105(1), pp.109-142.\\[-20pt]

\bibitem{Guivarch.Goldsheid2} Guivarc'h, Y. On the spectrum of a large subgroup of a semisimple group. Journal of modern dynamics, 2008, 2(1), p.15.\\[-20pt]

\bibitem{Guivarch.Raugi} Guivarc'h, Y. and Raugi, A. Fronti\`{e}re de Furstenberg, propri\'{e}t\'{e}s de contraction et th\'{e}or\`{e}mes de convergence. Probability Theory and Related Fields, 1985, 69(2), pp.187-242.\\[-20pt]

\bibitem{LePage} Le Page, E. Th\'{e}or\`{e}mes limites pour les produits de matrices al\'{e}atoires. In Probability measures on groups. Springer Berlin Heidelberg, 1982, pp. 258-303.\\[-20pt]

\bibitem{Prasad} Prasad, G. $\mathbb{R}$-regular elements in Zariski dense subgroups. The Quarterly Journal of Mathematics, 1994, 45(4), pp.542-545.\\[-20pt]

\bibitem{Sert.tez} Sert C. Joint spectrum and large deviation principle for random matrix products, Phd. thesis, Universit\'{e} Paris-Sud, 2016. \\[-20pt]

\bibitem{Sert.LDP.Zariski} Sert C. Large deviation principle for random matrix products (in preparation).\\[-20pt]

\bibitem{Sert.growth.indicator} Sert C. Growth indicator on semisimple linear groups (in preparation).\\[-20pt]

\bibitem{Sert.joint.spectrum} Sert C. Joint spectrum on semisimple linear groups (in preparation).\\[-20pt]

\bibitem{Sert.hyperbolic.LDP} Sert C. Large deviation principle and growth indicator for Gromov hyperbolic groups (in preparation).\\[-20pt]

\bibitem{Tutubalin} Tutubalin, V.N. A central limit theorem for products of random matrices and some of its applications. In Symposia Mathematica Vol. 21, 1977, pp.101-116.
\end{thebibliography}
\end{document}